\documentclass[12pt]{amsart}
\usepackage{amscd,amssymb}
\usepackage[graph,frame,poly,arc]{xy}  
\usepackage{tikz}
\usepackage[plainpages,backref,urlcolor=blue]{hyperref}

\theoremstyle{plain}
\newtheorem{thm}[subsection]{Theorem}
\newtheorem{lem}[subsection]{Lemma}
\newtheorem{prop}[subsection]{Proposition}
\newtheorem{cor}[subsection]{Corollary}

\theoremstyle{definition}
\newtheorem{rk}[subsection]{Remark}

\newtheorem{ex}[subsection]{Example}

\numberwithin{equation}{section}
\newcommand{\HH}{{\mathcal H}}

\newcommand{\A}{{\mathcal A}}
\newcommand{\B}{{\mathcal B}}

\newcommand{\D}{{\mathcal D}}
\newcommand{\CC}{{\mathcal C}}

\newcommand{\al}{{\alpha}}

\newcommand{\C}{\mathbb{C}}
\newcommand{\PP}{\mathbb{P}}

\DeclareMathOperator{\rank}{rank}

\begin{document}

\title [Waring rank and singularities of some projective hypersurfaces]
{Waring rank of symmetric tensors, and singularities of some projective hypersurfaces}

\author[Alexandru Dimca]{Alexandru Dimca$^1$}
\address{Universit\'e C\^ ote d'Azur, CNRS, LJAD and INRIA, France and Simion Stoilow Institute of Mathematics,
P.O. Box 1-764, RO-014700 Bucharest, Romania}
\email{dimca@unice.fr}

\author[Gabriel Sticlaru]{Gabriel Sticlaru}
\address{Faculty of Mathematics and Informatics,
Ovidius University
Bd. Mamaia 124, 900527 Constanta,
Romania}
\email{gabrielsticlaru@yahoo.com }

\thanks{$^1$ This work has been partially supported by the French government, through the $\rm UCA^{\rm JEDI}$ Investments in the Future project managed by the National Research Agency (ANR) with the reference number ANR-15-IDEX-01 and by the Romanian Ministry of Research and Innovation, CNCS - UEFISCDI, grant PN-III-P4-ID-PCE-2016-0030, within PNCDI III.}



\begin{abstract} 
We show that if a homogeneous polynomial $f$ in $n$ variables has Waring rank $n+1$, then the corresponding projective hypersurface $f=0$ has at most isolated singularities, and the type of these singularities is completely determined by the combinatorics of a hyperplane arrangement naturally associated with the Waring decomposition of $f$. We also discuss the relation between the Waring rank and the type of singularities on a plane curve, when this curve is defined by the suspension of a binary form, or when the Waring rank is 5.
\end{abstract}
 
\maketitle


\section{Introduction} 
For the general question of symmetric tensor decomposition we refer to
\cite{Guide, Car1+, Ci, CGLM, FOS, FLOS, IK,L1,LT,MO,O}, as well as to the extensive literature
quoted at the references in \cite{Guide} and \cite{FLOS}.
We describe first a possibly new general approach to tensor decompositions, and then illustrate this approach by a number of very  simple situations. Consider the graded polynomial ring $S=\C[x_1,...,x_n]$, let $f \in S_d$ be a homogeneous polynomial of degree $d$, such that the corresponding hypersurface

\begin{equation}
\label{eq1}
V=V(f):f=0
\end{equation}
in the complex projective space $\PP^{n-1}$ is reduced. We consider the {\it Waring decomposition}
\begin{equation}
\label{eq2}
(\D) \ \  \  \ \  f=\ell_1^d + \cdots +\ell_r^d,
\end{equation}
where $\ell_j \in S_1$ are linear forms, and $r$ is minimal, in other words $r=\rank f$ is the {\it Waring rank}  of $f$. We assume in the sequel that the form $f$ {\it essentially involves the $n$ variables $x_1, \ldots,x_n$}, in other words that $f$ cannot be expressed as a polynomial  in a fewer number of variables than $n$. This is equivalent to the fact that
the linear forms $\ell_j$'s span the vector space $S_1$, in particular we have $r \geq n$. When such a decomposition is given, we will also use the notation $V=V_{\D}$ to show that the hypersurface $V$ comes from the decomposition $(\D)$. For a given form $f\in S_d$ with Waring rank $r$, one can define the {\it Waring locus } of $f$ to be the set
\begin{equation}
\label{eq2.1}
W_f=\{ [\ell_1] \in \PP(S_1) \ : \ \text{ there are } \ell_2, \ldots \ell_r  \in S_1 \text{ such that } (1.2) \text{ holds } \}.
\end{equation}
 Moreover, the {\it locus of forbidden linear forms} is the complement $F_f=\PP(S_1) \setminus W_f$, see \cite{BalChi,Car1+}.
Consider the linear embedding
\begin{equation}
\label{eq3}
\varphi _{\D}:\PP^{n-1} \to \PP^{r-1}, \ \ x \mapsto ( \ell_1(x):  \cdots :\ell_r(x)),
\end{equation}
determined by the decomposition $(\D)$. In the projective space $\PP^{r-1}$ we have two basic objects, namely the {\it Fermat hypersurface} of degree $d$, given by
\begin{equation}
\label{eq4}
F:f_F(y)=y_1^d+ \cdots +y_r^d=0,
\end{equation}
and the {\it Boolean arrangement}
\begin{equation}
\label{eq5}
\B:f_{\B}(y)=y_1y_2 \cdots y_r=0.
\end{equation}
With these notations, we clearly have
\begin{equation}
\label{eq6}
V_{\D}=\varphi_{\D}^{-1}(F).
\end{equation}
Alternatively, let $E_{\D}$ be the $(n-1)$-dimensional linear subspace in $\PP^{r-1}$ given by the image of $\varphi _{\D}$, and note that we have an isomorphism
\begin{equation}
\label{eq7}
V_{\D} \simeq F \cap E_{\D}.
\end{equation}
Hence, to understand the geometry of the hypersurface $V_{\D}$ in terms of the tensor decomposition $(\D)$, we have to analyze the position of the linear subspace $E_{\D}$ with respect to the Fermat hypersurface $F$. One way to do this, is to consider the hyperplane arrangement
\begin{equation}
\label{eq8}
\A_{\D}=\varphi_{\D}^{-1}(\B) : \ell_1  \cdots \ell_r=0
\end{equation}
in $\PP^{n-1}$, associated to the decomposition $(\D)$. The fact that the linear forms $\ell_j$'s span the vector space $S_1$ implies that $\A_{\D}$ is an {\it essential} arrangement, see \cite{DHA, OT} for general facts on hyperplane arrangements.
Since $\A_{\D}$ is nothing else but the intersection $\B \cap E_{\D}$, that is the trace of the arrangement $\B$ on the linear space
$E_{\D} \simeq \PP^{n-1}$, it follows that the position of $E_{\D}$ is reflected in the properties of this induced arrangement $\A_{\D}$. Our general idea is to fix the combinatorics of the hyperplane arrangement $\A_{\D}$, e.g. by fixing the intersection lattice of the corresponding central arrangement, and see which geometric properties of the hypersurface $V_{\D}$ can be derived just from this combinatorics. 
Note that a similar idea, namely the study of the Fano scheme $F_k(X_{r,d})$ of projective $k$-planes contained in the projective hypersurface in $\PP^{rd-1}$ given by
$$X_{r,d}: \sum_{i=1}^r\prod_{j=1}^dx_{ij}=0,$$
was used already by N. Ilten, H. S\"u\ss{ }and Z.Teitler, see \cite{IT,IS}, to study the decompositions of a homogeneous polynomial $f$ as a sum of products of linear forms. The equation of the hypersurface $X_{r,d}$ can be regarded as a polarization of the equation \eqref{eq4} for the Fermat hypersurface.

In this paper we illustrate this approach by  three simple and hopefully interesting cases. The first one is when $r=n+1$. Indeed, recall that by our assumption $r \geq n$, and note that the case $r=n$ is rather trivial, i.e. in this case the hypersurface $V_{\D}$ is projectively equivalent to the Fermat hypersurface $F$ of degree $d$. The main result in this case is Theorem \ref{thm1H}, saying that the 
hypersurface $V_{\D}$ has at most isolated singularities, and the type of the corresponding singularities is determined by the combinatorics of the hyperplane arrangement $\A_{\D}$. The fact that the singularities of the hypersurface  $V_{\D}$ are at most isolated follows also from a very general result due to Landsberg and
Teitler, see \cite[Theorem 1.3]{LT}.
On the other hand, the number of these singularities is not determined by the combinatorics of the hyperplane arrangement $\A_{\D}$, but by the geometry of the hypersurface $\hat F$, which is the dual of the Fermat hypersurface $F$.
More precisely, when $r=n+1$, then
$E_{\D}$ is a hyperplane in $ \PP^n$,  the hypersurface
$V_{\D}$ is singular exactly when $E_{\D} \in \hat F$, and the number of singularities of $V_{\D}$ is equal to the number of local irreducible components of $\hat F$ at the corresponding point $E_{\D}$. 
Among the nodal hypersurfaces constructed in this way are the 
{\it generalized Cayley hypersurfaces} discussed in Example \ref{ex1H},
with additional information for {\it generalized Cayley curves} in
Proposition \ref{prop1H}. We end the discussion of the forms $f$ of Waring rank  $r=n+1$ with a result relating the Waring locus $W_f$ of the form $f$, or rather its complement $F_f$, to the dual variety of the hypersurface $V$, see Proposition \ref{propWR}.

The second case is when $n=3$, the Waring rank $r$ is arbitrary, but the line arrangement $\A_{\D}$ has the simplest combinatorics, i.e. $\A_{\D}$ has a point of multiplicity $r-1$. The main result in this case is Corollary \ref{corS}, which shows again that the possible singularities of $V_{\D}$
in this situation are very restricted.

\begin{figure}[h]
\centering
\begin{tikzpicture}[scale=0.7]
\draw[style=thick,color=blue] (-0.5,4) -- (1,-1);
\draw[style=thick,color=blue]  (-1,-1) -- (0.5,4);
\draw[style=thick,color=blue] (-2,0) -- (2,0);
\draw[style=thick,color=blue] (0,4) -- (0,-1);
\draw[style=thick,color=blue] (-2,-1) -- (1,4);
\node at (0,-1.5) {$n_4=1, n_2=4$};
\node at (0,-2.5) {$(1)$};
\end{tikzpicture}
\hspace*{0.1in}
\begin{tikzpicture}[scale=0.7]
\draw[style=thick,color=blue] (-0.5,4) -- (1,-1);
\draw[style=thick,color=blue]  (-1,-1) -- (0.5,4);
\draw[style=thick,color=blue] (-2,0) -- (2,0);
\draw[style=thick,color=blue] (0,4) -- (0,-1);
\draw[style=thick,color=blue] (-2,-1) -- (2,1);
\node at (0,-1.5) {$n_3=2, n_2=4$};
\node at (0,-2.5) {$(2)$};
\end{tikzpicture}
\hspace*{0.1in}
\begin{tikzpicture}[scale=0.7]
\draw[style=thick,color=blue] (-0.5,4) -- (1,-1);
\draw[style=thick,color=blue]  (-1,-1) -- (0.5,4);
\draw[style=thick,color=blue] (-2,0) -- (2,0);
\draw[style=thick,color=blue] (-2,4) -- (2,-1);
\draw[style=thick,color=blue] (0,4) -- (0,-1);
\node at (0,-1.5) {$n_3=1, n_2=7$};
\node at (0,-2.5) {$(3)$};
\end{tikzpicture}
\hspace*{0.1in}
\begin{tikzpicture}[scale=0.7]
\draw[style=thick,color=blue] (-0.5,4) -- (1,-1);
\draw[style=thick,color=blue]  (-1,-1) -- (0.5,4);
\draw[style=thick,color=blue] (-2,0) -- (2,0);
\draw[style=thick,color=blue] (-2,4) -- (2,-1);
\draw[style=thick,color=blue] (-2,-1) -- (2,4);
\node at (0,-1.5) {$n_2=10$};
\node at (0,-2.5) {$(4)$};
\end{tikzpicture}


\caption{Five lines in the plane}
\label{fig:5lines}
\end{figure}
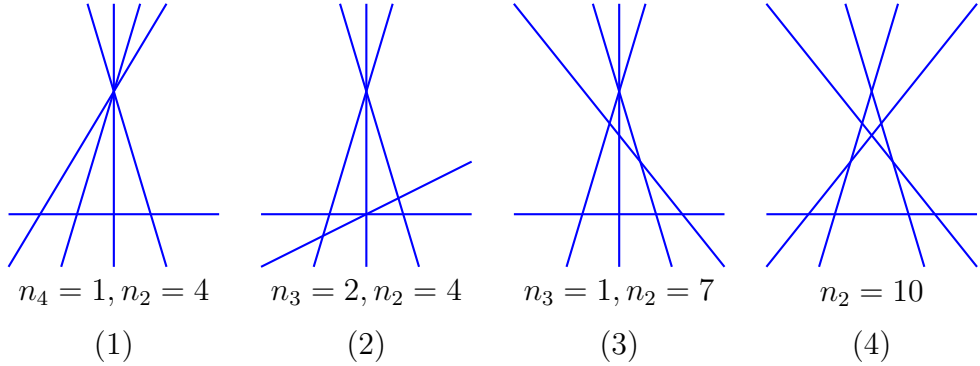

Finally we consider the plane curves of Waring rank 5. In this case the combinatorics of the line arrangement $\A_{\D}$ displays four possibilities, as shown in Figure \ref{fig:5lines}. In the first two cases, our results are complete, see Corollary \ref{corS} and Proposition \ref{propT2}, while in the other two cases we can for the moment give only partial results, see Proposition \ref{propT3} and Proposition \ref{propT4}.

The first author thanks AROMATH team at INRIA Sophia-Antipolis for excellent working conditions, and  Laurent Bus\'e, Bernard Mourrain and Alessandro Oneto for stimulating discussions. 

We also thank Zach Teitler for drawing our attention to some very interesting related results and references, and
Joerg Fricke for the one-line proof of Lemma \ref{lemK}.
Finally, we thank
the referee, for the very careful reading of the manuscript and for the very useful suggestions and remarks.

\section{ The hypersurfaces in $\PP^{n-1}$ of Waring rank $r=n+1$, with $n \geq 3$}

Let $\mu(V)$ (resp. $\tau(V)$) denote the global Milnor (resp. global Tjurina number) of a projective hypersurface $V$ having only isolated singularities, that is the sum of all the local Milnor numbers $\mu(V,p)$ (resp. local Tjurina numbers $\tau(V,p)$) over all the singular points $p \in V$. For the definition of local Milnor number and Tjurina number, we refer to \cite{DRCS}.

Consider the Fermat hypersurface $F$ defined in \eqref{eq4}, as well as its {\it dual hypersurface} $\hat F$ in $\hat \PP^{r-1}$. If $z_1,...,z_r$ denote the coordinates on $\hat \PP^{r-1}$, then the defining equation 
$$\hat F: \Delta_F(z)=0,$$
of the dual hypersurface $\hat F$ can be obtained by eliminating
$y_1,...,y_r$ in the following system of equations
$$z_j=y_j^{d-1}  \text{ for } j=1,...,r \text{ and  } f_F(y)=0,$$
e.g. using the computer algebra software SINGULAR, see \cite{Sing}.
Note also that 
$$\deg \Delta_F=d(d-1)^{r-2},$$
see for instance \cite{Kl}.
When $r=n+1$, then
$E_{\D}$ is a hyperplane in $ \PP^n$, and the hypersurface
$V_{\D}$ is singular exactly when $E_{\D} \in \hat F$.  An example of this dual variety $\hat F$ is given below in Example \ref{ex3B}. 

We start with the following result on essential hyperplane arrangements.

\begin{prop}
\label{propHA} 
Let $\CC$ be a central, essential  hyperplane arrangement in $\C^n$,
consisting of $n+1$ hyperplanes. Then there is a unique integer $k$, satisfying $2 \leq k \leq n$, and such that $\CC$ is linearly equivalent to the hyperplane arrangement
$$\CC_k:x_1 \cdots x_n (x_1+ \cdots +x_k)=0.$$
\end{prop}
\proof
The arrangement $\CC$ being essential, we can choose $n$ hyperplanes
$H_1, \dots, H_n$ in $\CC$ such that $H_1 \cap \dots \cap H_n=0$.
If $\CC_0$ is the arrangement formed by these $n$ hyperplanes, then its intersection lattice $L(\CC_0)$ is isomorphic to the intersection lattice of the Boolean arrangement in $\C^n$. If $H$ is the hyperplane in $\CC$ distinct from the $H_j$, $j=1, \dots ,n$, then let $k$ the smallest codimension of a flat $X \in L(\CC_0)$ which is contained in $H$.
We can assume that $X=H_1 \cap \dots \cap H_k$. It remains to choose the linear coordinates on $\C^n$ such that $H_j:x_j=0$ for $j=1, \dots ,n$. Then the equation of $H$ has the form
\begin{equation}
\label{eqH}
H: a_1x_1+ \cdots +a_kx_k=0,
\end{equation}
where all $a_j$ are non-zero, by the choice of $k$.  The claim follows, replacing $x_j$ by $a_jx_j$ for $j=1, \cdots, k$.
\endproof

Next we prove a technical results, needed in the proof of Theorem \ref{thm1H}.

\begin{lem}
\label{lemH1}
Let 
$$f=x_1^d+ \cdots +x_n^d+ (a_1x_1+\dots + a_kx_k)^d,$$
where  $2 \leq k \leq n$ and $a=(a_1, \dots, a_k) \in (\C^*)^k$. Assume that either $n \geq 4$ or $n=3$ and $d \geq 4$. Then the Hessian polynomial
$hess(f)$ is divisible by  a form $\al^{d-1}$, with  $\al \in S_1$, if and only if $\al$ is proportional to one of the linear forms
 $x_j$ for $j=k+1, \ldots, n$.
\end{lem}
\proof
It is clear that the Hessian polynomial
$hess(f)$ is divisible by $x_j^{d-2}$ for $j=k+1, \ldots, n$. To prove the remaining claim, it is enough to show that, for $k=n$, the Hessian polynomial
$hess(f)$ is not divisible by any form $\al^{d-1}$, with $\al \in S_1$.
It is easy to see that, when $k=n$, one has the following formula
\begin{equation}
\label{eqHess}
hess(f)(x)=d^n (d-1)^n\left(\prod_{j=1}^nx_j^{d-2}\right)\left( 1+\ell^{d-2} \left(\sum_{j=1}^n \frac{a_j^2}{x_j^{d-2}}\right) \right),
\end{equation}
where $\ell=\sum_{j=1}^n a_jx_j$. This formula shows that  the Hessian polynomial
$hess(f)$ is  divisible neither by some $x_i$, nor by the linear form $\ell$.
Suppose that $hess(f)$ is  divisible by some linear form $\al$, not proportional to any $x_i$ or to  $\ell$. We can assume by symmetry that $x_1$ occurs with a non-zero coefficient in $\al$, and so we can write
$$\al(x)=x_1-\sum_{j=2}^nb_jx_j,$$
with at least one coefficient $b_j \ne 0$. Let $H_0:\al=0$.
If $\al$ divides $hess(f)$, then it follows that the restriction of $hess(f)$ to the hyperplane $H_0$ is trivial, that is
$$hess(f)(\sum_{j=2}^nb_jx_j,x_2,\ldots,x_n)=0.$$
Using the formula \eqref{eqHess}, we get
\begin{equation}
\label{eqHess2}
1+(\ell(\sum_{j=2}^nb_jx_j,x_2,\ldots,x_n))^{d-2} \left(\frac{a_1^2}{(\sum_{j=2}^nb_jx_j)^{d-2}}+\sum_{j=2}^n \frac{a_j^2}{x_j^{d-2}}\right)=0.
\end{equation}
Identify the hyperplane $H_0$  to $\C^{n-1}$, with coordinates $x_2,\ldots,x_n$.
Consider next the hyperplanes in this $\C^{n-1}$, given by
$$H: \ell(\sum_{j=2}^nb_jx_j,x_2,\ldots,x_n)=0,  \  \  H_1: \sum_{j=2}^nb_jx_j=0$$
and
$H_j:x_j=0$ for $j=2, \ldots,n$. Note that these last $n-1$ hyperplanes are distinct.

\medskip {\bf Case 1: $n \geq 4$.} Then at least one of the $n-1$ hyperplanes $H_2,...,H_n$ is distinct from both $H_1$ and $H$. Let's say $H_n \ne H_1$ and $H_n \ne H$. Then, choose a point $p \in H_n \setminus (H_1 \cup H_2 \cup ... \cup H_{n-1} \cup H)$ and note that, by taking the limit when 
$$ x'=(x_2,...,x_n) \to p, \  \ x' \notin (H_1 \cup H_2 \cup ... \cup H_{n-1} \cup H_n \cup H)$$
we get a contradiction using the equation \eqref{eqHess2}.

\medskip {\bf Case 2: $n=3$, $d \geq 4$.} The new situation that may occur here is when $H$ is one of the two lines $H_2$ and $H_3$, and $H_1$ is the other line in the set $\{H_2,H_3\}$. By symmetry we can assume that $H_1=H_2$ and $H=H_3$. The equality $H_1=H_2$ implies $b_2 \ne 0$ and $b_3=0$. The equality $H=H_3$ implies $b_2=-a_2/a_1$. If we replace these data in \eqref{eqHess2} we get
$$1+(a_3x_3)^{d-2}\left( \left(\frac{a_1^2}{b_2^{d-2}}+ a_2^2\right)\frac{1}{x_2^{d-2}}+\frac{a_3^2}{x_3^{d-2}}\right)=0.$$
If the coefficient
$$c=\frac{a_1^2}{b_2^{d-2}}+ a_2^2$$
is non-zero, we get a contradiction as above by taking $x'=(x_2,x_3) \to (0,1)$, with $x_3 \ne 0$. Consider now the case $c=0$. Then one gets
$$a_3^d+1=0.$$
Using all these relations, a direct computation shows that
$$\frac{hess(f)(x)}{d^3 (d-1)^3}=\left( \frac{a_2^2}{a_1^{d-2}}\right)x_3^{d-2}\ell^{d-2}\left( (a_1x_1)^{d-2}-(-a_2x_2)^{d-2}\right)+$$
$$+a_3^2(x_1x_2)^{d-2}\sum_{j=0,d-3}{d-2 \choose j} a_3^jx_3^j(a_1x_1+a_2x_2)^{d-2-j}.$$
This formula implies that $hess(f)(x)$ is divisible by 
$$\al=a_1x_1+a_2x_2=a_1(x_1-b_2x_2),$$
 but not by $\al ^2$, which completes the proof of our claim.
\endproof
\begin{ex}
\label{ex10H}
Note that for
$$f=-3(x+y)(y+z)(x+z)=x^3+y^3+z^3-(x+y+z)^3$$
one has $hess(f)=-6^3(x+y)(y+z)(x+z)$. This shows that the case $n=d=3$ is indeed exceptional and has to be excluded in Lemma \ref{lemH1}.
\end{ex}
Here is the main result of this paper.

\begin{thm}
\label{thm1H}
Consider the reduced hypersurface $V=V_{\D}:f=0$ of degree $d$ in $\PP^{n-1}$, of Waring rank $n+1$, with $n \geq 3$.
Then, up-to a linear change of coordinates, there is a unique integer $k$,
satisfying $2 \leq k \leq n$, and such that
$$V_{\D}: f=x_1^d+ \cdots +x_n^d+ (a_1x_1+\dots + a_kx_k)^d=0,$$
where $a=(a_1, \dots, a_k) \in T^k=(\C^*)^k$. Moreover, 
the following hold.
\begin{enumerate}

\item The projective hypersurface $V_{\D}$ is singular if and only if $R_k(a)=0$, where $R_k(a)$ is the resultant of the system of $k$ equations $(S)$, with $(k-1)$ indeterminates $u=(u_1,...,u_{k-1})$, given by
$$g_j(u)=a_ku_j^{d-1}-a_j=0 \text{ for } j=1,...,k-1$$
and
$$g_{k}(u)=a_k(a_1u_1+...+a_{k-1}u_{k-1}+a_k)^{d-1}+1=0.$$

\item When the hypersurface $V_{\D}$ is singular, then it has only isolated singularities of type 
$A_{2^{k-1},d^{n-k}}$, given in local coordinates by the equation
$$v_1^2+ ...+v_{k-1}^2+w_1^d+...+w_{n-k}^d=0.$$ 
In particular, the hypersurface $V_{\D}$ is irreducible for $n \geq 4$.
\item These singularities are located at the points $(u^0:1:0: ...:0)\in \PP^{n-1}$, where $u^0=(u_1^0,...,u_{k-1}^0)$ is a solution of the system of equations $(S)$. In particular,
$$\mu(V_{\D})=\tau(V_{\D})=N(S)(d-1)^{n-k},$$
 where $N(S)$ is the number of solutions of the  system of equations $(S)$.

\end{enumerate}

\end{thm}

\proof By assumption, the polynomial $f$ has a Waring decomposition $(\D)$ as in \eqref{eq2},
with $r=n+1$. Then the associated hyperplane arrangement $\A_{\D}$,
or the central version of it in $\C^n$ to be more precise, satisfies the conditions in Proposition \ref{propHA}. Then formula \eqref{eqH} implies that the defining equation $f=0$ can be chosen as claimed.
The unicity of the integer $k$ follows from Lemma \ref{lemH1}. Indeed, if there is another linear coordinate system $y=(y_1,...,y_n)$ on $\C^n$ such that
$$f(y)=y_1^d+ \cdots + y_nd+(b_1y_1+ \cdots +b_my_m)^d,$$
with $b=(b_1, ..., b_m) \in (\C^*)m$, then passing from the coordinates $x$ to the coordinates $y$ just multiplies the Hessian polynomial of $f$ by a constant, so the divisibility properties by forms of the type $\al^{d-2}$ are preserved. Lemma \ref{lemH1} tells us that $k=m$.
When $V:f=0$ is singular, then the  unicity of $k$ follows straightforward  from the claim (2).

The proof of the claim (1) is by direct computation,  using the system of equations given by the vanishing of all first order partial derivatives of $f$.
Any solution $x^0=(x^0_1,...,x^0_n)$ of this system satisfies
$x^0_j=0$ if and only if $j>k$. We set $x_k=1$ and $u_j=x_j$ for $1 \leq j <k$, and get in this way the system (S).

We give two proofs for the claim (2): the first one is by a direct but lengthy elementary computation, while the second one is geometrical, using the properties of the inflection points of the Fermat hypersurface $F$.

\medskip

{\bf The first proof of the claim (2)}
Fix a solution $u^0=(u_1^0,...,u_{k-1}^0)$ of the system (S) and consider the corresponding singular point of $V_{\D}$, namely $p=(u^0:1:0: ...:0)$. We choose the local coordinates $(v,w)$ at
$p$ such that
$x_j=u^0_j+v_j$ for $1 \leq j <k$, $x_k=1$ and  $x_m=w_{m-k}$ for $k<m\leq n$. Then the singularity $(V,p)$ is given by the germ at the origin of the  polynomial $h(v,w)=h_1(v)+h_2(w)$ in $(v,w)$, where
$$h_1(v)=(u^0_1+v_1)^d+...+(u_{k-1}^0+v_{k-1})^d+1+(\gamma+a_1v_1+...+a_{k-1}v_{k-1})^d$$
where $\gamma=a_1u^0_1+...+a_{k-1}u_{k-1}^0+a_k$,
and $h_2(w)=w_1^d+...+w_{n-k}^d.$ It is easy to check that the polynomial $h_1(v)$ has only terms of degree $\geq 2$. To complete the proof, it is enough to show that the quadratic form $q=j^2h_1$, given by degree two part in $h_1$, is non-degenerate, see if necessary the first pages in \cite{DRCS}.
If we omit the binomial coefficient ${d \choose 2}$ which is a common factor, the quadratic form $q$ is given by
$$q=(u_1^0)^{d-2}v_1^2+...+(u_{k-1}^0)^{d-2}v_{k-1}^2+\gamma^{d-2}
(a_1v_1+...+a_{k-1}v_{k-1})^2.$$
The system (S) implies 
$$(u_j^0)^{d-2}=\frac{a_j}{a_ku_j^0}$$
pour $1 \leq j <k$ and also
$$\gamma^{d-2}=-\frac{1}{a_k\gamma}.$$
Note that the system implies that $\gamma \ne 0$, so all denominators are non zero. We have to show that the $(k-1)\times (k-1)$ symmetric matrix $M(q)$ associated to the quadratic form $q$ has a non zero determinant. If we multiply all the elements in this matrix by $\gamma a_k$, we get a new matrix $N(q)$ with elements
$n_{i,j}=-a_ia_j$ if $i \ne j$ and 
$$n_{j,j}=\frac{a_j\gamma}{u_j^0}-a_j^2,$$
for $1 \leq j <k$.
Multiply the $j$-th row in this matrix by $u_j^0$, for $1 \leq j <k$, and call the resulting rows $L_1,...,L_{k-1}$. Next add all the rows $L_j$ with $1<j<k$ to the first row $L_1$, and get in this way the row
$$L_1'=(a_1a_k,a_2a_k,...,a_{k-1}a_k).$$
Then divide by $a_k$ and get the new first row
$$L''_1=(a_1,a_2,...,a_{k-1}).$$
If we add $a_ju_j^0L''_1$ to the row $L_j$, we get a new row $L'_j$, where all the elements are zero, except the diagonal element which is
$a_j\gamma \ne0$. Hence we have shown by this sequence of elementary transformations on the rows, that the matrix $M(q)$ is non degenerate. This ends the first proof of the claim (2).

\medskip

{\bf The second proof of the claim (2)} It is clearly enough to consider the case $k=n$.
 It is known that the set of inflection points of any hypersurface $V$ is given by the intersection of $V$ and its Hessian hypersurface $\HH_V$. For the Fermat hypersurface, the Hessian hypersurface, with reduced structure, is given by
$$\HH_F: y_1y_2 \cdots y_{n+1}=0.$$
Now choose a point $p \in F$. If all the coordinates of $p$ are non zero, then $p$ is not an inflection point, and hence there are two possibilities for a plane $H$ passing through $p$. Either $H=T_pF$, and then the 
hypersurface singularity $(F\cap H,p)$ in $(H,p)=(\C^{n-1},0)$ is an $A_1$-singularity, since $p$ is not an inflection point. Or else $H \ne T_pF$, and then $H$ is transversal to $F$ at the point $p$.
Assume now that some coordinates in $p$ are zero. Then the tangent space at $p$ will have an equation
$$T_pF:b'_1y_1+...+b'_{n+1}y_{n+1}=0,$$
with some of the coefficients $b'_j$ equal to zero. When $k=n$, the hyperplane $E_{\D}$ is the image of the map
$$\varphi(x)=(x_1:...:x_n:a_1x_1+...+a_nx_n),$$
and hence, it is given by the equation
$$a_1y_1+...+a_ny_n-y_{n+1}=0.$$
This implies that $E_{\D} \ne T_pF$, and then $E_{\D}$ is transversal to $F$ at the point $p$. This completes the second proof of claim (2).

\medskip

The claim (3) is obvious, since for an isolated singularity $(V_{\D},p)$ of type $A_{2^{k-1},d^{n-k}}$ as defined above, one clearly has
$$\mu(V_{\D},p)=\tau(V_{\D},p)=(d-1)^{n-k}.$$
Note also that any solution $u^0$ of the system (S) is a simple solution, i.e. a solution with multiplicity one.
\endproof

\begin{rk}
\label{rk0H}
When $n=2$ and $d\geq 3$, then it is obvious that a rank 3 binary form $f$ can be written as $f(x_1,x_2)=x_1^d+x_2^d+(a_1x_1+a_2x_2)^d$ for some $(a_1,a_2) \in (\C^*)^2$. Moreover, it is easy to check that such a binary form can have only linear factors of multiplicity $\leq 2$, e.g. by using Proposition \ref{propS}. As an example, note that
$$f(x_1,x_2)=x_1^4+x_2^4+(x_1+x_2)^4=2(x_1^2+x_1x_2+x_2^2)^2$$
has two double factors. For a different approach to binary forms of rank 3 see \cite{Tok}.
\end{rk}
The final part of the proof of Theorem \ref{thm1H} yields the following.
\begin{cor}
\label{corI2}
Any hyperplane $H:b_1y_1+...+b_{n+1}y_{n+1}=0$ in $\PP^n$ with $n \geq 2$, such that
$b_j \ne 0$ for all $j$, is either transversal to the Fermat hypersurface
$F$, or it is tangent to $F$ at a finite number of points, say $N(H)$, such that at each such point
$p \in F \cap H$, the hypersurface singularity $(F\cap H,p)$ in $(H,p)=(\C^{n-1},0)$ is an $A_1$-singularity.
\end{cor}

\begin{rk}
\label{rk1H} Any isolated hypersurface singularity may occur on a hyperplane section of a smooth projective hypersurface, see
  \cite [Proposition (11.6)]{DRCS}. It is rather surprising that the hyperpane sections of the Fermat hypersurface $F$ yield only singularities of  very limited number of types, i.e. the singularities $A_{2^{k-1},d^{n-k}}$ introduced above.
  In the setting of  Corollary \ref{corI2}, the hyperplane section $F \cap H$ is a nodal hypersurface, whose nodes are exactly the tangency points between $H$ and $F$. This gives
  the obvious inequality
 $$N(H) \leq N(d,n-1),$$
 where $N(e,m)$ denotes the maximal number of singularities a nodal hypersurface of degree $e$ in $\PP^m$ can have. Explicit upper bounds for $N(e,m)$ were given by Varchenko, see \cite{Var}.

\end{rk}

\begin{rk}
\label{rk2H}
Note that the global Milnor number $\mu(V_{\D})$ coincides with the multiplicity of the dual hypersurface $\hat F$ at the point $E_{\D}$, see \cite[Proposition (11.24)]{DRCS}. Note also that the hypersurface $V_{\D}$ has only nodes as singularities when $k=n$, and this says that the hypersurface germ $(\hat F,E_{\D})$ is a union of smooth  components, see for instance the equivalent properties in  \cite [(11.33)]{DRCS}. More precisely, to a node $p=(p_1:...:p_n) \in V_{\D}$, it corresponds the point 
$$q=\varphi_{\D}(p)=(p_1:...:p_n:a_1p_1+...a_np_n) \in F$$
such that $T_qF=E_{\D}$. Then the dual mapping $\phi: F\to \hat F$
sends the point $q$ to the point $\hat q=\phi(q)$ corresponding to the hyperplane $E_{\D}$.
And the corresponding smooth component $(Z,\hat q)$ of the hypersurface germ $(\hat F,\hat q)$ has a (projective) tangent space in $\PP^n$ given by
$$T_{\hat q}Z :p_1z_1+p_2z_2+...+p_nz_n+(a_1p_1+...a_np_n)z_{n+1}=0,$$
see the proof of \cite[Proposition (11.24)]{DRCS}.

\end{rk}

\begin{rk}
\label{rk3H}
With the notation from Theorem \ref{thm1H}, it follows that for any
$2 \leq k \leq n$, there is a positive integer $m_k>0$ such that
$$R_k(a_1,...,a_k)^{m_k}=\Delta_F(a_1,...,a_k,0,...,0,-1),$$
up-to a non-zero constant factor.
To see this, it is enough to notice that the hyperplane $E_{\D}$
corresponds to the point
$$(a_1:...:a_k:0:...:0:-1) \in \hat \PP^{n}.$$
See Example \ref{ex3B} for situations where $m_k=1$ and $m_k=2$.
\end{rk}

\begin{ex}[Generalized Cayley Hypersurfaces]
\label{ex1H}
Consider the reduced hypersurface $V_{\D}:f=0$ of odd degree $d$ in $\PP^{n-1}$, with $n \geq 3$, given by
$$V_{\D}: f=(n-2)^{d-1}(x_1^d+ \cdots +x_n^d)-(x_1+\dots + x_n)^d=0.$$
Then $V_{\D}$ has $n$ singularities $A_1$ located at the points
$p^{i}=(p^{i}_1, \cdots , p^{i}_n)$ for $i=1,...,n$, with $p^{i}_j=1$ for $i \ne j$ and $p^{i}_i=-1$. First note that $d$ odd implies $f_{x_j}(p^{i})=0$ for 
$i,j=1,...,n$, where $f_{x_j}$ denotes the partial derivative of $f$ with respect to $x_j$. Hence all the points $p^{i}$ are singular points of 
the hypersurface $V_{\D}$, and the fact that they are nodes $A_1$ follows from Theorem  \ref{thm1H}, case $k=n$. The classical Cayley surface corresponds to $n=4$ and $d=3$, see \cite{Cay}. The case $n=d=3$, when $V_{\D}$ is a triangle, is also discussed below in Example \ref{ex3B} (3).
\end{ex}
It is a challenging problem to describe {\it all} the singularities of a generalized Cayley hypersurface. We give the result only for plane curves below, see Proposition \ref{prop1H}.

\medskip

Finally in this section, we say something about the Waring locus $W_f$ as defined in \eqref{eq2.1} for forms $f$ of rank $n+1$ using Theorem \ref{thm1H}. The fact that the complement $F_f=\PP(S_1) \setminus W_f$ is related to the dual $\hat V$ has been noticed already in a number of cases, see for instance \cite[Theorem 3.16]{Car1+} for plane cubic curves. To state our result, we define
$\hat V'$ to be the set of point $\hat q \in \hat V$ for which there is at least an irreducible component of the analytic germ $(\hat V, \hat q)$ which is not a smooth germ.
\begin{prop}
\label{propWR}
Consider the reduced hypersurface $V:f=0$ of degree $d$ in $\PP^{n-1}$, of Waring rank $n+1$, with $n \geq 3$, given by
$$ f=x_1^d+ \cdots +x_n^d+ (a_1x_1+\dots + a_nx_n)^d=0,$$
where $a=(a_1, \dots, a_n) \in (\C^*)^n$.
Then $\hat V' \subset F_f$.
\end{prop}
\proof
Using the unicity of the integer $k$ in Theorem \ref{thm1H}, it follows that any Waring decomposition of $f$ has the form
$$ \D: \  \  \  f=\ell_1^d + \cdots +\ell_r^d$$
with $r=n+1$ and the linear forms $\ell_j$'s playing a symmetric role. In particular, there is a vector $b=(b_2, \dots, b_r) \in (\C^*)^n$ such that
$$\ell_1=b_2\ell_2+ \ldots + b_r \ell_r.$$
Let $H_1: \ell_1=0$ and note that the corresponding hyperplane section
$V \cap H_1$ can be identified, using the isomorphism \eqref{eq7}, to the intersection
$$F \cap E_{\D}\cap \{y_1=0\}.$$
The intersection $F\cap \{y_1=0\}$ is just the Fermat hypersurface
$$y_2^d+ \ldots + y_r^d=0$$
in the hyperplane $ \{y_1=0\}=\PP^{n-1}$, while $E_{\D}\cap \{y_1=0\}$ is the hyperplane
$$b_2y_2+ \ldots + b_r y_r=0$$
in $ \{y_1=0\}=\PP^{n-1}$. Using Corollary \ref{corI2}, it follows that 
the hyperplane section
$V \cap H_1$ is either smooth, i.e. $H_1 \notin \hat V$, or it has at most finitely many nodes $A_1$ as singularities. 
This yields our claim $\hat V' \subset F_f$, since the points in $\hat V'$
correspond precisely to hyperplanes $H$ in $\PP^{n-1}$ such that
$V \cap H$ has at least one singular point which is not a node, recall our discussion above in Remark \ref{rk2H}.
\endproof

\section{The case of plane curves of Waring rank $r=4$}
 
In this section we consider the case $n=3$ in more detail. 
We set $x_1=x$, $x_2=y$, $x_3=z$,  $a_1=a$, $a_2=b$ ,$a_3=c$,
$z_1=A$, $z_2=B$, $z_3=C$ and $z_4=D$
to simplify the notation. 
 The corresponding line arrangement $\A_{\D}$ in $\PP^2$ consists of 4 lines, not all of them passing through one point. It follows that there are two possibilities for the combinatorics of $\A_{\D}$: either $\A_{\D}$ has a triple point and 3 nodes, which is the case $k=2$ in Theorem \ref{thm1H}, or $\A_{\D}$ is a generic arrangement, and has 6 nodes,
which is the case $k=3$ in Theorem \ref{thm1H}.

\begin{figure}[h]
\centering
\begin{tikzpicture}[scale=0.9]
\draw[style=thick,color=blue] (-0.5,4) -- (1,-1);
\draw[style=thick,color=blue]  (-1,-1) -- (0.5,4);
\draw[style=thick,color=blue] (-2,0) -- (2,0);
\draw[style=thick,color=blue] (0,4) -- (0,-1);
\node at (0,-1.5) {$n_3=1, n_2=3$};
\node at (0,-2.5) {$(1)$};
\end{tikzpicture}
\hspace*{0.9in}
\begin{tikzpicture}[scale=0.9]
\draw[style=thick,color=blue] (-0.5,4) -- (1,-1);
\draw[style=thick,color=blue]  (-1,-1) -- (0.5,4);
\draw[style=thick,color=blue] (-2,0) -- (2,0);
\draw[style=thick,color=blue] (-2,4) -- (2,-1);
\node at (0,-1.5) {$n_3=0, n_2=6$};
\node at (0,-2.5) {$(2)$};
\end{tikzpicture}


\caption{Four lines in the plane}
\label{fig:4lines}
\end{figure}
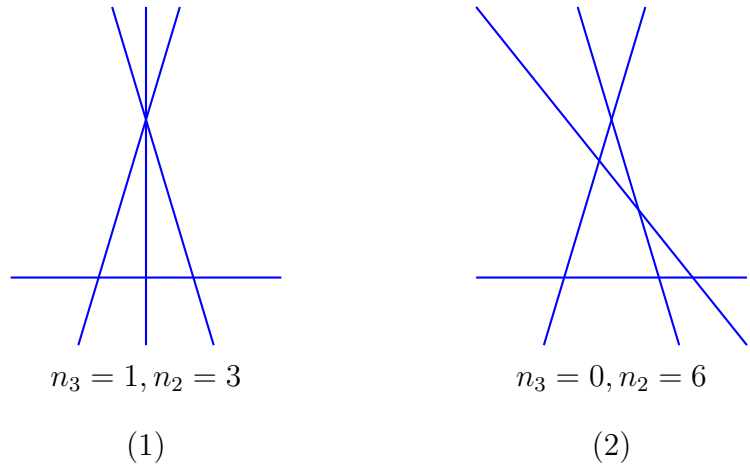

Then Theorem \ref{thm1H} yields  the following result.
\begin{cor}
\label{corI}
If the reduced plane curve $V_{\D}:f=0$ of degree $d$ has Waring rank $r=4$, then either $V_{\D}$ is smooth, or one of the following two cases occurs.
\begin{enumerate}

\item The line arrangement $\A_{\D}$ has a triple point and 3 nodes,
 the curve $V_{\D}$  has only simple singularities of type $A_{d-1}$. Moreover, $V_{\D}$ is irreducible for $d$ odd, and can have at most two irreducible components for $d$ even.

\item The line arrangement $\A_{\D}$ is generic, i.e. it has only nodes $A_1$, and then
 the curve $V_{\D}$  has also only simple singularities of type $A_1$.

\end{enumerate}

\end{cor}
Recall that a simple singularity of type $A_{d-1}$ is a singularity isomorphic to the singularity given by $v^2+w^d=0$, in the local coordinates $(v,w)$ at the origin of $\C^2$, see \cite{DRCS}. In particular, $A_1$ is a node, $A_2$ is a simple cusp, $A_3$ is a tacnode and $A_4$ is a ramphoid cusp.
With this explanations, the only claim in Corollary \ref{corI} that needs a proof is the claim about the number of irreducible components of $V_{\D}$ in case (1).
This follows from the following result.
\begin{thm}
\label{thm2A}
Consider the plane curve $$V_{\D}:f=x^d+y^d+z^d+(ax+by)^d,$$
where $(a,b) \in T^2=(\C^*)^2$. Then
the following hold.
\begin{enumerate}

\item The curve $V_{\D}$ is irreducible when $d$ is odd. 

\item When $d=2d'$ is even, then 
the curve $V_{\D}$ has $e \leq 2$ irreducible components. If $e=2$, then both components are smooth and 
$N(g_1,g_2) = d'.$
Moreover, when  $d'=2$, the equality $N(g_1,g_2) = d'$ implies that the curve $V_{\D}$ has $e=2$ components.
\end{enumerate}

\end{thm}

\proof

Recall that a plane curve $C:f=0$ has $e$ irreducible components if and only if $H^1(\PP^2 \setminus C,\C)$ is a $(e-1)$ dimensional vector space, see \cite[Proposition 4.1.13]{DSTH}. Moreover, one has
$$H^1(\PP^2 \setminus C,\C)=H^1(F_f,\C)_1,$$
where $F_f: f(x,y,z)=1$ is the Milnor fiber of $f$, and $H^1(F_f,\C)_1$ denotes the fixed part under the monodromy action, see \cite{DSTH} for details if necessary. Note that our polynomial $f$ can be written as
$$f(x,y,z)=h(x,y)+z^d,$$
where $h(x,y)=x^d+y^d+(ax+by)^d$. Now we use the formula \cite[(6.2.25)]{DSTH}, and conclude that $\dim H^1(\PP^2 \setminus C,\C)=e-1$ holds if and only if $\dim H^0(F_h,\C)=e$, where
$F_h:h(x,y)=1$ is the Milnor fiber of $h$. Then we use \cite[Proposition 3.2.3]{DSTH} and conclude that $\dim H^0(F_h,\C)=e$ if and only if
$h=h_1^e$, where $h_1 \in \C[x,y]$ is not the power of another polynomial. Note that, one has
$$bf_x-af_y=d(bx^{d-1}-ay^{d-1}),$$
so a polynomial with $(d-1)$ distinct roots.
On the other hand, one has
$$bf_x-af_y=bh_x-ah_y=e(bh_{1x}-ah_{1y})h_1^{e-1}.$$
This implies that $e \leq 2$.
When $d$ is odd, then only the case $e=1$ is possible, and hence $V_{\D}$ is irreducible in this case. Notice that, for $d$ odd, the singularity $A_{d-1}$ is unibranch, and hence, in particular, this gives another proof that the curve $V_{\D}$ is irreducible in this case.

Assume now that $d=2d'$ is even and $e=2$. Then
$$f(x,y,z)=h_1(x,y)^2+z^d=(h_1(x,y)+iz^{d'})(h_1(x,y)-iz^{d'}).$$
It follows that the curve $V_{\D}$ has two irreducible components, namely
$$C_1:f_1=h_1(x,y)+iz^{d'}=0 \text{ and } C_2:f_1=h_1(x,y)-iz^{d'}=0.$$
The two components intersect exactly at the points given by
$$h_1(x,y)=z=0.$$
At these points the two curves must be smooth, since the local singularities $A_{d-1}$ have two smooth branches, with a contact of order $d'$. It follows that $h_1(x,y)$ has only distinct roots, and that the number of singular points of $V_{\D}$ is exactly $d'$. This implies
$N(g_1,g_2) = d'=d/2$ in this case.
It remains to show that, conversely, when $N(g_1,g_2) = d'=d/2 =2$, then $V_{\D}$ is not irreducible.
Note that  the $\delta$-invariant of an $A_{2d'-1}$ singularity is 
$$\delta(A_{2d'-1})= \frac{(2d'-1)+2-1}{2}=d'.$$
If $V_{\D}$ is irreducible, then this would imply
$$N(g_1,g_2)\delta(A_{2d'-1})=(d' )^2 \leq \frac{(d-1)(d-2)}{2}=(d'-1)(2d'-1)$$
This inequality is impossible for $d'=2$. The case $d'=3$ is discussed in Example \ref{ex6A} below, where we show that $N(g_1,g_2) < d'$ for any choice of $(a,b)\in T^2$.
\endproof
 
First we discuss some examples in the case (1) of Corollary \ref{corI},
and hence we assume $k=2$ and $f=x^d+y^d+z^d+(ax+by)^d$ with $(a,b)\in T^2$.

\begin{ex}
\label{ex3A}
When $d=3$, the resultant $R_2(a,b)$ is the determinant of the following matrix
\begin{center}
$$M(g_1,g_2)=\left(
  \begin{array}{ccccccc}
     b & 0& -a & 0 \\
     0 &b& 0& -a  \\
    a^2b& 2ab^2 & b^3+1 & 0 \\
   0 & a^2b& 2ab^2 & b^3+1 \\
  
    \end{array}
\right).$$
\end{center}
Hence $V_{\D}$ is singular in this case if and only if
$$R_2(a,b)=\det M(g_1,g_2)=(a^3-b^3)^2+2(a^3+b^3)+1=0,$$
and in this case the curve has some cusps $A_2$. Since a cubic can have at most one cusp, it follows that the polynomials $g_1$ and $g_2$ have at most one root in common. This is reflected by the fact that the zero set in $\C^2$ of the ideal $I_3(M(g_1,g_2))$ generated by all the $3 \times 3$ minors of the matrix $M(g_1,g_2)$ is disjoint from the Zariski open set $T^2$.

\end{ex}

\begin{ex}
\label{ex4A}
When $d=4$, the resultant $R_2(a,b)$ is the determinant of the following matrix
\begin{center}
$$M(g_1,g_2)=\left(
  \begin{array}{ccccccc}
     b & 0& 0&-a & 0 &0\\
     0 &b& 0& 0& -a  &0 \\
0&0 &b& 0& 0& -a\\
    a^3b& 3a^2b^2 & 3ab^3& b^4+1 & 0&0 \\
   0 &  a^3b& 3a^2b^2 & 3ab^3& b^4+1 & 0 \\
  0&0& a^3b& 3a^2b^2 & 3ab^3& b^4+1 \\
    \end{array}
\right).$$
\end{center}
Hence $V_{\D}$ is singular in this case if and only if
$$R_2(a,b)=\det M(g_1,g_2)=(a^4+b^4)^3+3(a^8-7a^4b^4+b^8)+3(a^4+b^4)+1=0,$$
and in this case the curve has some singularities $A_3$. 
The set of pairs $(a,b)$ such that  the polynomials $g_1$ and $g_2$ have at least two common roots
is given by the zero set of the 
 ideal $I_5(M(g_1,g_2))$ generated by all the $5 \times 5$ minors of the matrix $M(g_1,g_2)$. Using the software SINGULAR \cite{Sing}, we see that this set has several irreducible components
which intersect the Zariski open set $T^2$, namely the points
$$(1,\pm 1), (\pm 1, \pm i), (\pm i, \pm i),$$
and the points obtained from these points using the transposition $(a,b) \mapsto (b,a)$. Here and in the sequel $i$ denotes a complex number with $i^2=-1$.
For any of these special values, the curve $V_{\D}$ has 2 singularities of type $A_3$.
 Note that the polynomial $g_1(t)$ has only simple roots for any degree $d$, hence the common roots of $g_1$ and $g_2$ are all distinct.
Moreover, the zero set in $\C^2$ of the ideal $I_4(M(g_1,g_2)$ generated by all the $4 \times 4$ minors of the matrix $M(g_1,g_2)$ is disjoint from the Zariski open set $T^2$, which is in accord with the obvious fact the a quartic cannot have more than two singularities $A_3$.
Note that the $\delta$-invariant of an $A_3$ singularity is $2$,  and hence a quartic curve with two $A_3$ singularities is reducible. It is easy to see that such a curve is the union of two smooth conics, tangent to each other in two points, corresponding to the two $A_3$ singularities. For instance, when $a=b=1$, we get
$$f=x^4+y^4+(x+y)^4+z^4=[\sqrt 2(x^2+y^2+xy)-iz^2][\sqrt 2(x^2+y^2+xy)+iz^2].$$
The corresponding smooth conics are tangent at the two points 
$(1:\al:0)$ with $\al^3=1$, $\al \ne 1$.
\end{ex}

\begin{ex}
\label{ex6A}
When $d=6$, the resultant $R_2(a,b)$ is the determinant of the following matrix $M(g_1,g_2)$ equal to
\begin{center}
$$\left(
  \begin{array}{cccccccccc}
b & 0 & 0 &0 & 0 & -a & 0 & 0 & 0 & 0 \\
0 & b & 0 & 0 & 0 & 0 & -a & 0 & 0 &0 \\
0&0&b&0&0&0&0&-a&0&0 \\
0&0&0& b&0&0&0&0&-a&0 \\
0&0&0&0& b& 0& 0& 0& 0 &-a\\
a^5b& 5a^4b^2& 10a^3b^3 & 10a^2b^4 & 5ab^5 & b^6+1&0&0&0&0\\
0 & a^5b& 5a^4b^2& 10a^3b^3 & 10a^2b^4 & 5ab^5 & b^6+1&0&0&0\\
0&0& a^5b& 5a^4b^2& 10a^3b^3 & 10a^2b^4 & 5ab^5 & b^6+1&0&0\\
0&0& 0& a^5b& 5a^4b^2& 10a^3b^3 & 10a^2b^4 & 5ab^5 & b^6+1&0\\
0&0&0&0&a^5b& 5a^4b^2& 10a^3b^3 & 10a^2b^4 & 5ab^5 & b^6+1\\
    \end{array}
\right).$$
\end{center}
Hence $V_{\D}$ is singular in this case if and only if
$$R_2(a,b)=\det M(g_1,g_2)=(a^6+b^6)^5+5(a^{24}-121a^{18}b^6+381a^{12}b^{12}-121a^6b^{18}+b^{24})+$$
$$+5(2a^{18}+381a^{12}b^6+381a^6b^{12}+2b^{18})+5(2a^{12}-121a^6b^6+2b^{12})
+5(a^6+b^6)+1=0,$$
and in this case the curve has some singularities $A_5$. 
The set of pairs $(a,b)$ such that  the polynomials $g_1$ and $g_2$ have at least two common roots
is given by the zero set of the 
 ideal $I_9(M(g_1,g_2))$ generated by all the $9 \times 9$ minors of the matrix $M(g_1,g_2)$. Using the software SINGULAR \cite{Sing}, we see that this set has several irreducible components
which intersect the Zariski open set $T^2$, for instance the points
$(a,1)$, where $a$ is a solution of the equation
$$a^{12}-11a^6-1=0,$$
are in this intersection.
On the other hand, the set of pairs $(a,b)$ such that  the polynomials $g_1$ and $g_2$ have at least three common roots
is given by the zero set of the 
 ideal $I_8(M(g_1,g_2))$ generated by all the $8 \times 8$ minors of the matrix $M(g_1,g_2)$. Using the software SINGULAR \cite{Sing}, we see that this set does not intersect the Zariski open set $T^2$.
 Hence the equality $N(g_1,g_2) = d'=d/2$ can not hold for $d' = 3$, and hence all the curves $V_D$ are irreducible in this case.
\end{ex}

Next we discuss some examples in the case (2) of Corollary \ref{corI},
and hence we assume $k=3$ and $f=x^d+y^d+z^d+(ax+by+cz)^d$, with $(a,b,c) \in T^3$.

\begin{ex}
\label{ex3B}
When $d=3$, the resultant $R_3(a,b,c)$, obtained using the command $elim$ in SINGULAR, is given by
$$R_3(a,b,c)=a^{12}-4a^9b^3+6a^6b^6-4a^3b^9+b^{12}-4a^9c^3+4a^6b^3c^3+4a^3b^6c^3-4b^9c^3+6a^6c^6+$$
$$+4a^3b^3c^6+6b^6c^6-4a^3c^9-4b^3c^9+c^{12}+4a^9-
4a^6b^3-4a^3b^6+4b^9-4a^6c^3+40a^3b^3c^3-$$
$$-4b^6c^3-4a^3c^6-4b^3c^6+4c^9+6a^6+4a^3b^3+6b^6+4a^3c^3+4b^3c^3+6c^6+4a^3+4b^3+4c^3+1.$$
Hence $V_{\D}$ is singular in this situation if and only if
$R_3(a,b,c)=0$, and in this case the curve $V_{\D}$ has a number of nodes $A_1$. The dual variety $\hat F$ of the Fermat hypersurface in $\PP^3$ is given by the equation
\begin{equation}
\label{eqD1}
\Delta(A,B,C,D)=A^{12}-4A^9B^3+6A^6B^6-4A^3B^9+B^{12}-4A^9C^3+4A^6B^3C^3+
\end{equation}
$$+4A^3B^6C^3-4B^9C^3+6A^6C^6+4A^3B^3C^6+6B^6C^6-4A^3C^9-4B^3C^9+C^{12}-$$$$-4A^9D^3+4A^6B^3D^3+4A^3B^6D^3-4B^9D^3+4A^6C^3D^3-40A^3B^3C^3D^3+4B^6C^3D^3+$$$$+4A^3C^6D^3+4B^3C^6D^3-4C^9D^3+6A^6D^6+4A^3B^3D^6+6B^6D^6+4A^3C^3D^6+4B^3C^3D^6+$$$$+6C^6D^6-4A^3D^9-4B^3D^9-4C^3D^9+D^{12}=0,$$
where $(A,B,C,D)$ are homogeneous coordinates on $\hat \PP^3$.
The plane $E_D$, which is the image of the map
$\varphi_D(x,y,z)=(x:y:z:ax+by+cz),$
is given by the equation 
$$ay_1+by_2+cy_3-y_4=0$$
in $\PP^3$, and hence corresponds to the point
$(A:B:C:D)=(a:b:c:-1).$
One can check that
\begin{equation}
\label{eqD2}
R_3(a,b,c)=\Delta(a,b,c,-1).
\end{equation}
As an example, consider the family of Waring decompositions $\D_a$
corresponding to the triple $(a,b,c)=(a,-a-2,-a-2)$. Then
$$R_3(a,-a-2,-a-2)=(a+1)^3(a^2-a+1)^2(25a^5+215a^4+841a^3+1777a^2+2015a+961).$$
The line $L \subset \hat \PP^3$ corresponding to the family of planes
$E_{\D_a}$ has the following intersection points with the dual hypersurface $\hat F$.
\begin{enumerate}

\item 5 simple points on $\hat F$, corresponding to the 5 roots of the irreducible factor $25a^5+215a^4+841a^3+1777a^2+2015a+961$.
For each such root $a$, the curve $V_{\D_a}$ is a nodal cubic.

\item 2 points of multiplicity two on $\hat F$, corresponding to the 2 roots of the irreducible factor $a^2-a+1$.
For each such root $a$, the curve $V_{\D_a}$ is the union of a smooth conic and a secant line. More precisely, one has in this case
$$f=(-6a+9)( y+z )$$$$(x^2+(2a-3)xy+(-18/7a+12/7)y^2+(2a-3)xz+(-38/7a+23/7)yz+(-18/7a+12/7)z^2).$$

\item one point of multiplicity three on $\hat F$, corresponding to the  root  of the irreducible factor $a+1$.
For  $a=-1$, the curve $V_{\D_{-1}}$ is a triangle. More precisely, one has in this case
$$f=-3(x+y)(x+z)(y+z)$$
and hence the 3 nodes are located at the points $p=(1:1:-1)$, $p'=(1:-1:1)$ and $p''=(-1:1:1)$. It follows, from the discussion in Remark \ref{rk2H}, that the 
hypersurface germ $(\hat F,\hat q)$ with $\hat q=E_{\D_{-1}}$ is a union of 3 smooth  components, say $(Z,\hat q)$, $(Z',\hat q)$ and $(Z'',\hat q)$, such that
$$T_{\hat q}Z: A+B-C-D=0, \ T_{\hat q}Z':A-B+C-D=0 \text{ and } T_{\hat q}Z'':-A+B+C-D=0.$$
Since for each root $a$ of $R_3(a,-a-2,-a-2)=0$ its multiplicity is equal to the global Milnor number $\mu(V_{\D_a})$, it follows that the line $L$ is
transverse to the smooth irreducible components of $\hat F$, at each point $\hat q=E_{\D_a}$, recall Remark \ref{rk2H}.
\end{enumerate}

Now we consider the relation between $R_2(a,b)$ and $R_3(a,b,c)$ when $d=3$. Note that
$$R_3(a,b,0)=a^{12}-4a^9b^3+6a^6b^6-4a^3b^9+b^{12}+
4a^9-
4a^6b^3-4a^3b^6+4b^9+$$
$$+6a^6+4a^3b^3+6b^6+4a^3+4b^3+1,$$
and hence one has $R_3(a,b,0)=R_2(a,b)^2$, where
$$R_2(a,b)=(a^3-b^3)^2+2(a^3+b^3)+1=0,$$
as in Example \ref{ex3A}. This fact can be explained as follows. The plane $E_D$, the image of the map
$\varphi_D(x,y,z)=(x:y:z:ax+by),$
is given by the equation 
$$ay_1+by_2-y_4=0$$
in $\PP^3$, and hence corresponds to the point
$(A:B:C:D)=(a:b:0:-1) \in \hat \PP^3.$
One can check that
\begin{equation}
\label{eqD3}
R_2(a,b)^2=\Delta(a,b,0,-1),
\end{equation} 
as in Remark \ref{rk3H}.
\end{ex}

\begin{rk}
\label{rk3B}
The Waring ranks for plane cubics are listed in \cite[Theorem 8.1]{LT} and in \cite{Car1+}, subsection (3.4). The only cubic curve of Waring rank $>4$ is the union of a smooth conic and a tangent line, where the Waring rank is 5. Note that this curve has a unique $A_3$ singularity, and is discussed in Example \ref{exT3} below.
\end{rk}

We end this section by describing all the singularities of the generalized Cayley plane curves.
\begin{prop}
\label{prop1H}
Consider the reduced plane curve $V_{\D}:f=0$ of odd degree $d\geq 3$ in $\PP^2$  given by
$V_{\D}: f=x^d+y^d+z^d-(x+y +z)^d=0.$

Then the curve $V_{\D}$ has exactly $3(d-2)$ singularities $A_1$ located at the points $p=(1:u:-u)$, $p'=(u:1:-u)$ and $p''=(u:-u:1)$, where $u^{d-1}=1$.

Moreover, for $d\geq 5$, the curve $V_{\D}$ has four irreducible components, namely the three lines $x+y=0$, $x+z=0$ and $y+z=0$, and a smooth curve of degree $d-3$ meeting each of these three lines in $d-3$ points, distinct from the vertices of the triangle $T:(x+y)(x+z)(y+z)=0$.
\end{prop}
\proof
Before starting the proof, we check that the number of nodes is correct.
Note that $u$ can take $(d-1)$ values, so at first sight we have
$3(d-1)$ singularities. But each of the points $(1:1:-1)$, $(1:-1:1)$ and $(-1:1:1)$ is in fact counted two times, e.g. $(1:1:-1)$ can be both $p$ and $p'$ for $u=1$.  Hence the total number of points is
$$3(d-1)-3=3(d-2),$$
since there are no other repetitions. 

To start the proof, note that $(x:y:z)$ is a singular point of $V_{\D}$ if and only if one has
$$x^{d-1}=y^{d-1}=z^{d-1}=(x+y+z)^{d-1}.$$
It follows that $x+y+z \ne 0$, and hence we normalize by setting
$$x+y+z=1.$$
The claim about the location of the singularities follows from Lemma \ref{lemK} below, which might be well known to specialists.
The fact that $x+y$ is a factor of $f$ follows by using the formula
$$(x+y+z)^d=\left((x+y)+z\right)^d$$
and the fact that $x^d+y^d$ is divisible by $x+y$, $d$ being odd.
By symmetry, it follows that
$$f=(x+y)(x+z)(y+z)g,$$
where $g$ is a homogeneous polynomial of degree $d-3$.
The reduced curve $g=0$ intersects the line $L:x+y=0$ in exactly $d-3$ simple points, since they should be nodes on the curve $V_{\D}$.
In this way we get $3(d-3)$ nodes for $V_{\D}$ situated on $g=0$,
in addition to the 3 nodes which are the vertices of the triangle $T$.
This shows that there no other singularities for the curve $V_{\D}$, and hence in particular the curve $g=0$ is smooth, and hence in particular irreducible.
\endproof
\begin{lem}
\label{lemK}
Let $u,v,w$ be three complex numbers on the unit circle such that
$$u+v+w=1.$$
Then at least one of them is equal to $1$.
\end{lem}
\proof
The following one-line proof was communicated to us by Joerg Fricke.
If we think about the complex numbers as vectors in the real plane, then
$0,u,u+v,u+v+w=1$ are the four vertices of a rhombus (which may be degenerated, i.e. all vertices on the real axis), because the four sides have the same length equal to 1.
\endproof
Note that this Lemma does not extend to four unitary complex numbers.
Indeed, if $u\ne 1$ satisfies $u^5=1$, then one clearly has
$$(-u)+(-u^2)+(-u^3)+(-u^4)=1.$$

\section{On the Waring rank of binary forms and the singularities of their suspensions}

In this section we recall first a property in the case $n=2$, i.e. when $f$ is a binary form in $x_1=x$ and $x_2=y$. A lot is known on the Waring rank of binary forms,  staring with the work of Sylvester \cite{Sy}, see also \cite{BM,CS,Re}, \cite[section 3.3]{Car1+}, \cite[Example 2.5]{MO}.
For the following useful result, due to Neriman Tokcan, see \cite[Theorem 3.1]{Tok}.

\begin{prop}
\label{propS}
Assume that the binary form $f$ has the following factorization
$$f=f_1^{m_1}\cdots f_s^{m_s},$$
for some $s \geq 2$, where the linear forms $f_i$ and $f_j$ are not proportional for $i \ne j$,
and $m_i\geq 1$. Let $m=\max \{m_i\}$.
Then the Waring rank of $f$ is at least $m+1$.
\end{prop}

Now we return to the case $n=3$.
Note that for a line arrangement $\A_{\D}$ coming from a Waring decomposition $(\D)$, the highest multiplicity of a point can be $r-1$, since the line arrangement is supposed to be essential. If there is such a point $p$, then the line arrangement consists just of $r-1$ line passing through $p$ and an additional secant line, i.e. the combinatorics of $\A_{\D}$ is the simplest possible. For Waring decompositions having this type of associated line arrangement we have the following result, which is a generalization of Corollary \ref{corI} (1) where $r=4$.
\begin{cor}
\label{corS}
If the reduced plane curve $V_{\D}:f=0$ of degree $d$ has Waring rank $r\geq 4$ and the line arrangement  $\A_{\D}$ has a  point of multiplicity $r-1$,
then 
 the curve $V_{\D}$  is either smooth or has only  singularities of type
 $v^m+w^d$, with $2 \leq m \leq r-2$. Moreover, in this case $r \leq d+1$.
\end{cor}
\proof
The hypothesis on the line arrangement  $\A_{\D}$ implies that the linear forms in the decomposition $(\D)$ can be chosen, up-to a linear change of coordinates, essentially as in the proof above, namely
$\ell_1=x$, $\ell_2=y$, $\ell_j=a_jx+b_jy$, for $j=3,...,r-1$, with
$a_jb_j\ne 0$ and the linear forms $\ell_i$ and $\ell_j$ are not proportional for $i \ne j$, and $\ell_r=z$. Then the binary form 
$$h(x,y)=\ell_1^d + \cdots +\ell_{r-1}^d$$
has only factors of multiplicity lower of equal to $r-2$ by Proposition \ref{propS}, and this proves our claim on the suspension $f=h(x,y)+z^d$.
The last claim follows from the fact that the Waring rank of a degree $d$ binary form is at most $d$, see \cite[Theorem 7.6]{Ge}.
\endproof
There is a similar result in higher dimensions, when the hyperplane arrangement $\A_{\D}$ has a codimension 2 edge $P$ which is the intersection of exactly $r-1$ hyperplanes in $\A_{\D}$, but we let the interested reader state this result for himself.

\section{On some plane curves of Waring rank 5}

An essential arrangement of 5 lines in the plane has one of the following pictures.
The plane curves of Waring rank 5 with the corresponding line arrangement $\A_{\D}$ of type (1) above were already discussed in Corollary \ref{corS}.
For the plane curves of Waring rank 5 with the corresponding line arrangement $\A_{\D}$ of type (2) above, we have the following result.
\begin{prop}
\label{propT2}
Let $V_{\D}:f=0$ be a reduced plane curve of degree $d$ and Waring rank $r=5$, such that the corresponding line arrangement $\A_{\D}$ is of type (2) above, i.e. there are $n_3=2$ triple points. Then, up-to a linear change of coordinates, we have
$$f=x^d+y^d+z^d+(a_1x+b_1y)^d+(a_2x+b_2z)^d,$$
with $(a_1,b_1,a_2,b_2) \in T^4=(\C^*)^4$.  Then the curve $V_{\D}:f=0$ is either smooth, or it has only nodes $A_1$ as singularities.
\end{prop}
\proof
The proof is by a direct computation, very similar to the first proof of Claim (2) in Theorem \ref{thm1H}, and is left to the reader.
Note that, in the special case $b_1^d+1=b_2^d+1=0$, there are solutions with $x=0$, and they should be treated in a separate way.
\endproof
For the plane curves of Waring rank 5 with the corresponding line arrangement $\A_{\D}$ of type (3) above, we show next that some new singularity types may occur.
\begin{prop}
\label{propT3}
Let $V_{\D}:f=0$ be a reduced plane curve of degree $d$ and Waring rank $r=5$, such that the corresponding line arrangement $\A_{\D}$ is of type (3) above, i.e. there are $n_3=1$ triple points. Then, up-to a linear change of coordinates, we have
$$f=x^d+y^d+z^d+(a_1x+b_1y)^d+(a_2x+b_2y+c_2z)^d,$$
with $(a_1,b_1,a_2,b_2,c_2) \in T^5=(\C^*)^5$ and $a_2b_1\ne a_1b_2$.  If the additional conditions
$$b_2^{d-1}+a_1(a_1b_2-a_2b_1)^{d-1}=(-1)^db_1b_2^{d-2}+a_1a_2^{d-2} =c_2^d+1=0 $$
hold, then the curve $V_{\D}:f=0$  has a singularity of type
$A_{2d-3}$ located at the point $p=(-b_2:a_2:0)$.
\end{prop} 
\proof
We look for conditions that our curve $V_{\D}$ has a singular point $p=(p_1:p_2:0)$ on the line $z=0$. Then the condition $f_z(p)=0$ implies that $a_2p_1+b_2p_2=0$, and hence $p=(-b_2:a_2:0)$ is the unique possibility.
The condition $f_x(p)=0$ yields
$$b_2^{d-1}+a_1(a_1b_2-a_2b_1)^{d-1}=0.$$
Similarly, the condition $b_1f_x(p)-a_1f_y(p)=0$ yields
$$(-1)^db_1b_2^{d-2}+a_1a_2^{d-2}=0.$$
Hence if these conditions are fulfilled, $p$ is a singularity of the curve
$V_{\D}:f=0$. Now we can write $p=(u:1:0)$ with $u=-b_2/a_2$ and use local coordinates at $p$ given by $x=u+v$, $y=1$ and $z=w$. Then the local equation of the germ $(V_{\D},p)$ is given by
$$f(u+v,1,w)=(u+v)^d+1+w^d+(a_1u+b_1+a_1v)^d+(a_2v+c_2w)^d.$$
If we expand this polynomial in $v,w$, the terms of degree $\leq 1$ vanish.
The coefficient $A$ of $v^2$ is given, after division by the coefficient ${d \choose 2}$,  by
$$u^{d-2}+a_1^2(a_1u+b_1)^{d-2}=-\frac{a_1}{u}(a_1u+b_1)^{d-1} +a_1^2(a_1u+b_1)^{d-2}=-\frac{a_1b_1}{u}(a_1u+b_1)^{d-2} \ne 0,$$
since $a_1u+b_1=-(a_1b_2-a_2b_1)/a_2$. Next we look at the terms of degree $d$. If the coefficient of $w^d$, which is $c_2^d+1$ is non-zero, then the singularity $(V_{\D},p)$ is of type $A_{d-1}$. But if this coefficient is zero, as assumed in our hypothesis, then $(V_{\D},p)$ is given by a semi-weighted homogeneous equation with leading term
$$Av^2+Bvw^{d-1},$$
with $B=a_2c_2^{d-1} \ne 0$, see \cite[Section (7.3)]{DRCS}. It follows that $(V_{\D},p)$ is analytically equivalent to the singularity 
$$A_{2d-3}: v^2+w^{2d-2}=0.$$
\endproof
\begin{ex}
\label{exT3}
Consider the cubic curve of Waring rank 5 given by
$$V_{\D}:f=x^3+y^3+z^3-\frac{1}{4}(x+y)^3+(x-y-z)^3=0.$$
This curve satisfies all the assumptions in Proposition \ref{propT3} and we have
$$f=1/4(x-y)(7x^2-8xy+y^2-12xz+12yz+12z^2).$$
It follows that $V_{\D}$ is the union of a smooth conic and a tangent line
at the point $(1:1:0)$, in accord with Remark \ref{rk3B}.
\end{ex}
For the plane curves of Waring rank 5 with the corresponding line arrangement $\A_{\D}$ of type (4) above, we have the following partial result.
\begin{prop}
\label{propT4}
Let $V_{\D}:f=0$ be a reduced plane curve of degree $d$ and Waring rank $r=5$, such that the corresponding line arrangement $\A_{\D}$ is of type (4) above, i.e. there are only double points. Then, up-to a linear change of coordinates, we have
$$f=x^d+y^d+z^d+(a_1x+b_1y+c_1z)^d+(a_2x+b_2y+c_2z)^d,$$
with $(a_1,b_1,c_1,a_2,b_2,c_2) \in T^6=(\C^*)^6$ and 
$$m_{12}=a_1b_2- a_2b_1 \ne 0,\  m_{13}=a_1c_2- a_2c_1 \ne 0, \  m_{23}=b_1c_2-b_2c_1\ne 0.$$
  Then the curve $V_{\D}:f=0$  can have only singularities $p=(p_1:p_2:p_3)$ of type
$A_{m}$. Moreover, if $p_1p_2p_3=0$, then $p$ can be only an $A_1$-singularity.
\end{prop} 

\proof
Any singular point $p=(p_1:p_2:p_3)$ satisfies the equation
$$ m_{23}p_1^{d-1}-m_{13}p_2^{d-1}+m_{12}p_3^{d-1}=0.$$
This implies that at most one coordinate $p_j$ can be zero.
If this is the case, then by symmetry we can choose $p_1=0$ and $p_3=1$, and a direct computation as in the first proof of Claim (2) in Theorem \ref{thm1H} gives the result in this case.
Assume now that $p_1p_2p_3\ne0$. By symmetry we can suppose $p_3=1$ and we can check by a direct computation that the vanishing of all the terms of degree $\leq 2$ in $u,v$ in the polynomial
$f(p_1+u,p_2+v,1)$ leads to a contradiction.
\endproof
\begin{rk}
\label{rkT4}
We do not know whether singularities $A_m$ with $m\geq 2$ can really occur in the setting of Proposition \ref{propT4}.
\end{rk}


\begin{thebibliography}{00}





\bibitem{BalChi} E. Ballico, L. Chiantini, Sets Computing the Symmetric Tensor Rank, Mediterranean Journal of Mathematics 10 (2013), 643--654.

\bibitem{Guide} A. Bernardi, E. Carlini, M. V. Catalisano, A. Gimigliano, A. Oneto, The Hitchhiker guide to: Secant Varieties and Tensor Decomposition, Mathematics 6(2018), 314, DOI  10.3390/math6120314.

\bibitem{BM} L. Brustenga i Moncus\' i, S. K. Masuti, On the Waring rank of binary forms: the binomial formula and a dihedral cover of rank two forms, arXiv:1901.08320.

\bibitem{Car1+} E. Carlini, M.V. Catalisano, A. Oneto,  Waring loci and the Strassen conjecture. Adv. Math. 314(2017), 630--662.

\bibitem{Cay} A. Cayley,  A Memoir on Cubic Surfaces. Philos. Trans. R. Soc. Lond., Ser. A 159 (1869), 231--326.


\bibitem{Ci} C. Ciliberto, Geometric aspects of polynomial interpolation in more variables and of Waring’s problem. In: European Congress of Mathematics, Barcelona 2000, pages 289--316. Springer, 2001.

\bibitem{CS} G. Comas and M. Seiguer, On the rank of a binary form. Foundations of Computational
Mathematics, 11(1)52011), 65--78.

\bibitem{CGLM} P. Comon, G. Golub, L.H. Lim, B. Mourrain, Symmetric tensors and symmetric tensor
rank,  SIAM Journal on Matrix Analysis and Applications, 30(2008),1254--1279.


\bibitem
{Sing} { W. Decker, G.-M. Greuel, G. Pfister \and H. Sch{\"o}nemann.} \newblock {\sc Singular} {4-0-1} --- {A} computer algebra system for polynomial computations, available at {http://www.singular.uni-kl.de} (2014).




\bibitem{DRCS}  A. Dimca,   {\em Topics on Real and Complex Singularities}, Vieweg Advanced Lecture in 
Mathematics, Friedr. Vieweg und Sohn, Braunschweig, 1987.

\bibitem{DSTH}  A. Dimca,   {\em Singularities and Topology of Hypersurfaces} , Universitext, Springer Verlag, New York,
1992.


\bibitem{DHA}  A. Dimca,   {\em Hyperplane Arrangements: An Introduction}, Universitext, Springer, 2017



















  













\bibitem{FOS} R. Fr\" oberg, G. Ottaviani,  B. Shapiro, On the Waring problem for polynomial rings Proceedings
of the National Academy of Sciences, 109 (2012), 5600--5602.

\bibitem{FLOS} R. Fr\" oberg,  S. Lundqvist, 
A. Oneto, 
B. Shapiro, 
Algebraic stories from one and from the other pockets, Arnold Math. J. 4 (2018),  137--160.

\bibitem{Ge} A.V. Geramita. Inverse systems of fat points: Waring’s problem, secant varieties of Veronese
varieties and parameter spaces for Gorenstein ideals. In: The Curves Seminar at Queen’s,
volume 10, pages 2--114, 1996.





\bibitem{IK} A. Iarrobino, V. Kanev, {\it Power Sums, Gorenstein Algebras, and Determinantal Loci}, Springer Lecture Notes 1721, 1999.

\bibitem{IT} N. Ilten, Z.Teitler,  Product ranks of the $3 \times 3$ determinant and permanent. Canad. Math. Bull. 59 (2016), 311--319.


\bibitem{IS} N. Ilten, H. S\"u\ss, Fano schemes for generic sums of products of linear forms, arXiv:1610.06770.



\bibitem{Kl} S.L. Kleiman,  The enumerative theory of singularities, in: Real and Complex
Singularities (Oslo 1976), Sijthoff and Noordhoff, Amsterdam 1977,
pp. 297-- 396.

\bibitem{L1} J.M. Landsberg, {\it Tensors: Geometry and Applications}, Graduate Studies in Mathematics vol.128, American Mathematical Soc., 2012.

\bibitem{LT} J.M. Landsberg, Z. Teitler, On the ranks and border ranks of symmetric tensors, Found. Comput. Math. 10(3) (2010) 339--366.



 



\bibitem{MO} B. Mourrain, A. Oneto,  On minimal decompositions of low rank symmetric tensors, arXiv:1805.11940.



\bibitem{O} A. Oneto, Waring type problems for polynomials, Doctoral Thesis in Mathematics at Stockholm University, Sweden, 2016.


\bibitem{OT} P. Orlik and H. Terao, {\em Arrangements of Hyperplanes,} Springer-Verlag, Berlin Heidelberg New York, 1992.


\bibitem{Re} B. Reznick,  On the length of binary forms. In: Quadratic and higher degree forms, Dev. Math., vol. 31, pp. 207--232. Springer, New York (2013).














\bibitem{Sy} J.J. Sylvester. Lx. on a remarkable discovery in the theory of canonical forms and of hyperdeterminants.
The London, Edinburgh, and Dublin Philosophical Magazine and Journal of
Science, 2(12)(1851), 391--410.

\bibitem{Tok} N. Tokcan, On the Waring rank of binary forms, 	arXiv:1610.09065,  to appear in Linear Algebra and Its Applications.

\bibitem{Var}  A.N. Varchenko, On the semicontinuity of the spectrum and an upper bound for the number of singular points of a projective hypersurface, J. Soviet Math. 270(1983), 735--739.




\end{thebibliography}
\end{document}